\documentclass[12pt,reqno]{amsart}
\usepackage{mathrsfs}
\usepackage{amsfonts}
\usepackage{latexsym,amsmath,amsfonts,amssymb,amsthm}
 \theoremstyle{plain}
\newtheorem{Thm}{Theorem}
\newtheorem{Lem}{Lemma}
\newtheorem{Cor}{Corollary}

\newtheorem*{Prob}{Problem}
\theoremstyle{definition}
\newtheorem*{Ack}{Acknowledgment}
\theoremstyle{remark}
\newtheorem*{Rem}{Remark}

\def\Z{\mathbb Z}
\def\N{\mathbb N}

\def\A{\mathcal{A}}
\def\sA{\mathscr{A}}
\def\P{\mathcal P}

\def\1{{\bf 1}}
\def\pmod #1{\ ({\rm mod}\ #1)}
\def\mod #1{\ {\rm mod}\ #1}
\def\ord {{\rm ord}}
\def\jacob #1#2{\genfrac{(}{)}{}{}{#1}{#2}}

\begin{document}

\title{On the integers of the form $p^2+b^2+2^n$ and $b_1^2+b_2^2+2^{n^2}$}
\author{Hao Pan}
\email{haopan79@yahoo.com.cn}
\author{Wei Zhang}
\email{zhangwei\_07@yahoo.com.cn}
\address{Department of Mathematics, Nanjing University,
Nanjing 210093, People's Republic of China} 
\maketitle
\begin{abstract}
We prove that the sumset
$$
\{p^2+b^2+2^n:\, p\text{ is prime and }b,n\in\N\}
$$
has a positive lower density. We also construct a residue class with
odd modulo, which contains no integer of the form $p^2+b^2+2^n$. And
similar results are established for the sumset
$$
\{b_1^2+b_2^2+2^{n^2}:\, b_1,b_2,n\in\N\}.
$$
\end{abstract}

\section{Introduction}
\setcounter{equation}{0} \setcounter{Thm}{0} \setcounter{Lem}{0}
\setcounter{Cor}{0}

Let $\P$ denote the set of all primes. In 1934, Romanoff
\cite{Romanoff34} proved that the sumset
$$
\mathcal{S}_1=\{p+2^n:\, p\in\P,\ n\in\N\}
$$
has a positive lower density. Subsequently van der Corput
\cite{vanderCorput50} proved the complement of $\mathcal{S}_1$, i.e.
$\N\setminus\mathcal{S}_1$, also has a positive lower density. In
fact, Erd\H os \cite{Erdos50} showed that every positive integer $n$
with $n\equiv 7629217\pmod{11184810}$ is not of the form $p+2^n$.
The key ingredient of Erd\H os' proof is to find a finite class of
residue classes with distinct moduli, which covers all integers.
Nowadays, Erd\H os' idea has been greatly extended, and for the
further related developments, the readers may refer to
\cite{Crocker71,CohenSelfridge75,Chen00,Sun00,Chen01a,Chen01b,SunLe01,SunYang02,Chen03,
Yuan04,Chen05,LucaStanica05,Chen07,ChenFengTemplier08,WuSun}.

In 1999, with help of Br\"udern and Fouvry's estimations on sums of
squares \cite{BrudernFouvry94}, Liu, Liu and Zhan
\cite{LiuLiuZhan99} proved a Romanoff-type result:

\medskip
{\it The sumset
$$
\mathcal{S}_2=\{p_1^2+p_2^2+2^{n_1}+2^{n_2}:\, p_1,p_2\in\P,\
n_1,n_2\in\N\}
$$
has a positive lower density.}
\medskip

\noindent The key of their proof is the following lemma:

\medskip
{\it For $1\leq m\leq N$,
$$
|\{p_1^2+p_2^2-p_3^2-p_4^2=m:\, p_i\in\P,\ p_i^2\leq
N\}\ll\mathfrak{S}_-(m)\frac{N}{(\log N)^4},
$$
where $\mathfrak{S}_-$ will introduced in Section 2.}
\medskip

\noindent In the other direction, recently Crocker \cite{Crocker08}
proved that there exist infinitely many positive integers not
representable as the sum of two squares and two (or fewer) powers of
$2$.

Motivated by all these results, in the present paper, we shall study
the sumset
$$
\mathcal{S}_3=\{p^2+b^2+2^{n}:\, p\in\P,\ n,n_2\in\N\}
$$
First, we have the following Romanoff-type result.
\begin{Thm}
\label{sumset} The set $\mathcal{S}_3$ has a positive lower density.
\end{Thm}
Next, we need to say something about the complement of
$\mathcal{S}_3$. It is not difficult to see that almost all integers
in $\mathcal{S}_3$ are of the form $4k+1$ or $8k+2$. However, we
shall prove that
\begin{Thm}
\label{sumsetcomplement} There exists a residue class with odd
modulo, which contains no integer of the form $p^2+b^2+2^n$.
\end{Thm}
Since the modulo in Theorem \ref{sumsetcomplement} is odd, by the
Chinese remainder theorem, clearly both the two sets
$$
\{x\in\N: x\equiv 1\pmod{4},\ x\not\in\mathcal{S}_3\}
$$
and
$$
\{x\in\N: x\equiv 2\pmod{8},\ x\not\in\mathcal{S}_3\}
$$
have positive lower densities.

Furthermore, we also have a similar result on the integers not of
the form $p^2+b^2-2^n$.
\begin{Thm}
\label{sumdiffcomplement} There exists a residue class with odd
modulo, which contains no integer of the form $p^2+b^2-2^n$.
\end{Thm}

A well-known result due to Landau \cite{Landau08} asserts that
$$
\{b_1^2+b_2^2\leq N:\, b_1,b_2\in\N\}=\frac{KN}{\sqrt{\log
N}}(1+o(1))
$$
where
$$
K=\frac{1}{\sqrt{2}}\prod_{\substack{p\in\P\\
p\equiv
3\pmod{4}}}\bigg(1-\frac{1}{p^2}\bigg)^{-\frac{1}{2}}=0.764223653\ldots.
$$
And obviously
$$
\{n\in\N:\, 2^{n^2}\leq N\}\ll\sqrt{\log N}.
$$
These facts suggests us to obtain the following results.
\begin{Thm}
\label{sumtwosqaures} The sumset
$$
\mathcal{S}_4=\{b_1^2+b_2^2+2^{n^2}:\, b_1,b_2,n\in\N\}
$$
has a positive lower density. And conversely there also exists a
residue class with odd modulo, which contains no integer of the form
$b_1^2+b_2^2+2^{n^2}$.
\end{Thm}

The proofs of Theorem \ref{sumset} and the first assertion of
Theorem \ref{sumtwosqaures} are applications of sieve method. And we
shall construct a suitable cover of $\Z$ with odd moduli to prove
Theorem \ref{sumsetcomplement}, Theorem \ref{sumdiffcomplement} and
the second assertion of Theorem \ref{sumtwosqaures} . Throughout our
proof, the implied constants by $\ll$, $\gg$ and $O(\cdot)$ will be
always absolute.

\section{Proof of Theorem \ref{sumset}}
\setcounter{equation}{0} \setcounter{Thm}{0} \setcounter{Lem}{0}
\setcounter{Cor}{0}

For ${\bf d}=(d_1,d_2,d_3,d_4)$ with $\mu({\bf
d}):=\mu(d_1)\mu(d_2)\mu(d_3)\mu(d_4)\not=0$, define
$$
A(m,q,{\bf d})=q^{-4}\sum_{\substack{1\leq a\leq q\\
(a,q)=1}}e(-am/q)S(q,ad_1^2)S(q,ad_4^2)S(q,-ad_2^2)S(q,-ad_3^2)
$$
and
$$
\mathfrak{S}(m,{\bf d})=\sum_{q=1}^\infty A(m,q,{\bf d}),
$$
where
$$
S(q,a)=\sum_{x=1}^qe(ax^2/q)
$$
and $e(\alpha)=\exp(2\pi\sqrt{-1}\alpha)$. In particular, we set
$\mathfrak{S}_-(m)=\mathfrak{S}(m,(1,1,1,1))$ and $\omega({\bf
d},m)=\mathfrak{S}(m,{\bf d})/\mathfrak{S}_-(m)$. By the arguments
in \cite[Eq. (8.7)]{LiuLiuZhan99}, we know
$$
\omega({\bf d},n)=\prod_{\substack{p^u\| d_1d_4\\ p^v\|
d_2d_3}}\omega_{u,v}(p),
$$
where $p^\beta\| m$ means $p^\beta\mid m$ but $p^{\beta+1}\nmid m$.
\begin{Lem}[Liu, Liu and Zhan {\cite[Lemma 8.1]{LiuLiuZhan99}}]
\label{llz} Suppose that $p\geq 3$, $p^u\| d_1d_2$ and $p^v\|
d_2d_3$. If $p\nmid m$, then
$$
\omega_{1,0}(p)=\begin{cases} p/(p-1),\quad\text{if
}\jacob{-m}{p}=1,\\
p/(p+1),\quad\text{if }\jacob{-m}{p}=-1,
\end{cases}
$$
$$
\omega_{0,1}(p)=\begin{cases} p/(p-1),\quad\text{if
}\jacob{m}{p}=1,\\
p/(p+1),\quad\text{if }\jacob{m}{p}=-1,
\end{cases}
$$
and $\omega_{1,1}(p)=p/(p+1)$. And if $p^\beta\|m$ for some
$\beta\geq 1$, then
$$
\omega_{1,0}(p)=\omega_{0,1}(p)=\frac{1+p^{-1}-p^{-\beta}-p^{-\beta-1}}{1+p^{-1}-p^{-\beta-1}-p^{-\beta-2}}
$$
and
$$
\omega_{1,1}(p)=\frac{3-p^{-1}-p^{1-\beta}-p^{-\beta}}{1+p^{-1}-p^{-\beta-1}-p^{-\beta-2}}.
$$
\end{Lem}

Let
$$
\sA=\{(x_1,x_2,x_3,x_4):\, x_1^2+x_4^2=x_2^2+x_3^2+m,\ 1\leq
x_i^2\leq N\},
$$
and for ${\bf d}=(d_1,d_2,d_3,d_4)$, let
$$
\sA_{{\bf d}}=\{(x_1,x_2,x_3,x_4)\in\sA:\, x_i\equiv0\pmod{d_i}\}.
$$

\begin{Lem}[Br\"udern and Fouvry {\cite[Theorem 3]{BrudernFouvry94}}, Liu, Liu and Zhan {\cite[Lemma 9.1]{LiuLiuZhan99}}]
\label{bfllz}
$$
|\sA_{{\bf d}}|=\frac{\omega({\bf
d},m)}{d_1d_2d_3d_4}\frac{\pi}{16}\mathfrak{S}_-(m)\mathfrak{I}(m/N)N+R(m,N,{\bf
d}),
$$
where
$$
\mathfrak{I}(\theta)=2\int_{\max\{0,-\theta\}}^{\min\{1,1-\theta\}}t^{-1/2}(1-\theta-t)^{1/2}dt
$$
and
$$
\sum_{\substack{d_1,d_2,d_3,d_4\leq D\\ {\bf d}=(d_1,d_2,d_3,d_4)\\
|\mu(d)|=1}}|R(m,N,{\bf d})|\ll N^{1-\epsilon}.
$$
\end{Lem}

Let $D=N^{1/30}$ and $z=N^{1/300}$. Define $$
P(z)=\prod_{\substack{p<z\\ p\text{ prime}}}p.
$$Let
$$
f(k)=|\{(x_1,x_2,x_3,x_4)\in\sA:\, x_1x_2=k\}|.
$$
\begin{Lem}
\label{selberg} For any $d\mid P(z)$ with $d\leq \sqrt{D}$,
\begin{align*}
\sum_{k\equiv
0\pmod{d}}f(k)=&\frac{\pi}{16}\mathfrak{S}_-(m)\mathfrak{I}(m/N)N\prod_{p\mid
d}\bigg(\frac{\omega_{1,0}(p)}{p}+\frac{\omega_{0,1}(p)}{p}-\frac{\omega_{1,1}(p)}{p^2}\bigg)\\&
+O\bigg(\sum_{\substack{d_1,d_2\mid d,\ d\mid d_1d_2\\
t_1\mid d/d_1,\ t_2\mid d/d_2\\ {\bf
d}=(d_1t_1,d_2t_2,1,1)}}|R(m,N,{\bf d})|\bigg).
\end{align*}
\end{Lem}
\begin{proof} Applying Lemma \ref{bfllz}, we have
\begin{align*}
\sum_{d\mid k}f(k)
=&\sum_{\substack{d_1,d_2\mid d\\
d\mid d_1d_2}}|\{(x_1,x_2,x_3,x_4)\in\sA:\, (x_1,d)=d_1,\
(x_2,d)=d_2\}|\\
=&\sum_{\substack{d_1,d_2\mid d\\
d\mid d_1d_2}}\sum_{\substack{(x_1,x_2,x_3,x_4)\in\sA\\
d_1\mid x_1,\ d_2\mid x_2}}\bigg(\sum_{t_1\mid
(x_1,d)/d_1}\mu(t_1)\bigg)\bigg(\sum_{t_2\mid (x_2,d)/d_2}\mu(t_2)\bigg)\\
=&\sum_{\substack{d_1,d_2\mid d\\
d\mid d_1d_2}}\sum_{\substack{t_1\mid d/d_1\\ t_2\mid d/d_2}}\mu(t_1)\mu(t_2)\sum_{\substack{(x_1,x_2,x_3,x_4)\in\sA\\
d_1t_1\mid x_1,\ d_2t_2\mid x_2}}1\\
=&\sum_{\substack{d_1,d_2\mid d,\ d\mid d_1d_2\\ t_1\mid d/d_1,\ t_2\mid d/d_2\\
{\bf d}=(d_1t_1,d_2t_2,1,1)}}\mu(t_1)\mu(t_2)\bigg(\frac{\omega({\bf
d},m)}{d_1t_1d_2t_2}\frac{\pi}{16}\mathfrak{S}_-(m)\mathfrak{I}(m/N)N+R(m,N,{\bf
d})\bigg).
\end{align*}
In view of Lemma \ref{llz},
\begin{align*}
&\sum_{\substack{d_1,d_2\mid d,\ d\mid d_1d_2\\ t_1\mid d/d_1,\ t_2\mid d/d_2\\
{\bf d}=(d_1t_1,d_2t_2,1,1)}}\mu(t_1)\mu(t_2)\frac{\omega({\bf
d},m)}{d_1t_1d_2t_2}\\=&\sum_{\substack{[d_1,d_2]=d\\ t_1\mid
d_2/(d_1,d_2)\\ t_2\mid d_1/(d_1,d_2)}}\mu(t_1)\mu(t_2)\prod_{p\mid
t_1t_2(d_1,d_2)}\frac{\omega_{1,1}(p)}{p^2}\prod_{p\mid
d_1/(d_1,d_2t_2)}\frac{\omega_{1,0}(p)}{p}\prod_{p\mid
d_2/(d_1t_1,d_2)}\frac{\omega_{0,1}(p)}{p}\\
=&\sum_{\substack{t_1t_2t_3t_4t_5=d}}\mu(t_1)\mu(t_2)\prod_{p\mid
t_1t_2t_3}\frac{\omega_{1,1}(p)}{p^2}\prod_{p\mid
t_4}\frac{\omega_{1,0}(p)}{p}\prod_{p\mid
t_5}\frac{\omega_{0,1}(p)}{p}\\
=&\prod_{p\mid
d}\bigg(\frac{\omega_{1,0}(p)}{p}+\frac{\omega_{0,1}(p)}{p}-\frac{\omega_{1,1}(p)}{p^2}\bigg).
\end{align*}
\end{proof}
Clearly
\begin{align*}
&|\{(x_1,x_2,x_3,x_4):\, x_1^2+x_4^2=x_2^2+x_3^2+m,\ 1\leq x_i^2\leq N,\ (x_1x_2,P(z))=1\}|\\
=&\sum_{(k,P(z))=1}f(k)\leq\sum_{k}f(k)\bigg(\sum_{d\mid
(k,P(z))}\lambda_d\bigg)^2=\sum_{d_1,d_2\mid
P(z)}\lambda_{d_1}\lambda_{d_2}\sum_{k\equiv0\mod{[d_1,d_2]}}f(k),
\end{align*}
where $\lambda_d$ are the weights appearing in Selberg's sieve
method with $\lambda_d=0$ for $d\geq z$ (cf. \cite[Chapter
3]{HalberstamRichert74}). In view of Lemma \ref{selberg},
\begin{align*}
&\sum_{d_1,d_2\mid P(z)}\lambda_{d_1}\lambda_{d_2}\sum_{k\equiv0\mod{[d_1,d_2]}}f(k)\\
=&\frac{\pi}{16}\mathfrak{S}_-(m)\mathfrak{I}(m/N)N\sum_{d_1,d_2\mid
P(z)}\lambda_{d_1}\lambda_{d_2}\prod_{p\mid
[d_1,d_2]}\bigg(\frac{\omega_{1,0}(p)}{p}+\frac{\omega_{0,1}(p)}{p}-\frac{\omega_{1,1}(p)}{p^2}\bigg)\\
&+\sum_{d_1,d_2\mid
P(z)}\lambda_{d_1}\lambda_{d_2}\cdot O\bigg(\sum_{\substack{d_1',d_2'\mid [d_1,d_2],\ [d_1,d_2]\mid d_1'd_2'\\
t_1\mid [d_1,d_2]/d_1',\ t_2\mid [d_1,d_2]/d_2'\\ {\bf
d}=(d_1t_1,d_2t_2,1,1)}}|R(m,N,{\bf d})|\bigg).
\end{align*}
By Selberg's sieve method, we have
$$
\sum_{d_1,d_2\mid P(z)}\lambda_{d_1}\lambda_{d_2}\prod_{p\mid
[d_1,d_2]}\bigg(\frac{\omega_{1,0}(p)}{p}+\frac{\omega_{0,1}(p)}{p}-\frac{\omega_{1,1}(p)}{p^2}\bigg)
=\frac{1}{G_1(z)},
$$
where
\begin{align*}G_1(z)=\sum_{\substack{d\mid P(z)\\
d<z}}\prod_{p\mid d}
\frac{\omega_{1,0}(p)p^{-1}+\omega_{0,1}(p)p^{-1}-\omega_{1,1}(p)p^{-2}}{1-\omega_{1,0}(p)p^{-1}-\omega_{0,1}(p)p^{-1}+\omega_{1,1}(p)p^{-2}}.
\end{align*}
And since $|\lambda_d|\leq 1$,
\begin{align*}
&\sum_{d_1,d_2\mid
P(z)}\lambda_{d_1}\lambda_{d_2}\cdot O\bigg(\sum_{\substack{d_1',d_2'\mid [d_1,d_2],\ [d_1,d_2]\mid d_1'd_2'\\
t_1\mid [d_1,d_2]/d_1',\ t_2\mid [d_1,d_2]/d_2'\\ {\bf
d}=(d_1t_1,d_2t_2,1,1)}}|R(m,N,{\bf d})|\bigg)\\
\ll&\sum_{\substack{d_1,d_2\mid P(z)\\ d_1,d_2<z^2\\
{\bf
d}=(d_1,d_2,1,1)}}\tau(d_1)^2\tau(d_2)^2\tau(d_1d_2)^2|R(m,N,{\bf
d})|\ll N^{1-\epsilon/2},
\end{align*}
where $\tau$ is the divisor function. Noting that
$\omega_{1,0}(p),\omega_{0,1}(p)=1+O(1/p)$ and
$\omega_{1,1}(p)=O(1)$, we have (cf. \cite[Lemma
4.1]{HalberstamRichert74})
$$
G_1(z)\gg\prod_{p<z}\bigg(1+\frac{\omega_{1,0}(p)}{p}+\frac{\omega_{0,1}(p)}{p}-\frac{\omega_{1,1}(p)}{p^2}\bigg)\gg(\log
z)^2.
$$
And it had been showed \cite[Eq. (2.9)]{LiuLiuZhan99} that
$$
\mathfrak{S}_-(m)\ll\prod_{\substack{p^\beta\| m\\ p\geq 3\\
\beta\geq
0}}\bigg(1+\frac{1}{p}-\frac{1}{p^{\beta+1}}-\frac{1}{p^{\beta+2}}\bigg).
$$
Finally,
\begin{align*}
&|\{(x_1,x_2,x_3,x_4)\in\sA:\, |x_1|<z\text{ or
}|x_2|<z\}|\\
\leq&\sum_{\substack{|x_1|<z,\ x_2^2\leq N\\\text{or }|x_2|<z,\
x_1^2\leq
N}}|\{(x_3,x_4):\,x_4^2-x_3^2=x_2^2-x_1^2+m\}|\\
\ll&\sum_{\substack{|x_1|<z,\ x_2^2\leq N\\\text{or }|x_2|<z,\
x_1^2\leq N}}\tau(x_2^2-x_1^2+m)\ll N^{2/3}z.
\end{align*}
Thus we obtain that
\begin{Thm}
\label{squaresumdiff} For a positive integer $m$, we have
\begin{align*}
&|\{(p_1,p_2,b_3,b_4):\, p_1^2+b_4^2=p_2^2+b_3^2+m,\ p_i\in\P,\
b_i\in\N,\ p_i^2,
b_i^2\leq N\}|\\
\ll&\frac{N}{(\log N)^2}\prod_{p\mid m}\bigg(1+\frac{1}{p}\bigg).
\end{align*}
\end{Thm}
Now we are ready to prove Theorem \ref{sumset}.
\begin{proof}[Proof of Theorem \ref{sumset}]
By the prime number theorem, clearly
$$
|\{p^2+b^2\leq N:\, p\in\P,\
b\in\N\}|\leq|\P\cap[0,\sqrt{N}]|\cdot|\N\cap[0,\sqrt{N}]|\ll\frac{N}{\log
N}.
$$
On the other hand, in \cite{Rieger68}, Rieger proved that
$$
|\{p^2+b^2\leq N:\, p\in\P,\ b\in\N\}|\gg\frac{N}{\log N}.
$$
Define
$$
r(x)=|\{(p,b,n):\,p^2+b^2+2^n=x,\ p\in\P,\ b,n\in\N\}|.
$$
Recall that $\mathcal{S}_3=\{x\in\N:\,r(x)\geq 1\}$. Then by the
Cauchy-Schwarz inequality,
$$
|\{(p,b,n):\,p^2+b^2+2^n\leq N,\ p\in\P,\ b,n\in\N\}|=\sum_{x\leq
N}r(x)\leq\sqrt{|\mathcal{S}_3\cap[1,N]|}\cdot\sqrt{\sum_{x\leq
N}r(x)^2}.
$$
Clearly
\begin{align*}
&|\{(p,b,n):\,p^2+b^2+2^n\leq N,\ p\in\P,\ b,n\in\N\}|\\
\geq&|\{p\in\P:\,p^2\leq N/3\}|\cdot|\{b\in\N:\,b^2\leq
N/3\}|\cdot|\{n\in\N:\,2^n\leq N/3\}|\gg N.
\end{align*}
So it suffices to
show that
$$
\sum_{x\leq N}r(x)^2\ll N.
$$
Applying Theorem \ref{squaresumdiff}, we have
\begin{align*}
\sum_{x\leq N}r(x)^2=&
|\{(p_1,p_2,b_1,b_2,n_1,n_2):\,p_1^2+b_1^2+2^{n_1}=p_2^2+b_2^2+2^{n_2}\leq
N\}|\\
\leq&2\sum_{\substack{2^{n_1}\leq 2^{n_2}\leq N}}
|\{(p_1,p_2,b_1,b_2):\,p_1^2+b_1^2-p_2^2-b_2^2=2^{n_2}-2^{n_1},\
p_i^2, b_i^2\leq N\}|\\
\ll&\frac{N}{\log N}\cdot\frac{\log N}{\log
2}+\sum_{2^{n_1}<2^{n_2}\leq N}\frac{N}{(\log N)^2}\prod_{p\mid
2^{n_2}-2^{n_1}}\bigg(1+\frac{1}{p}\bigg).
\end{align*}
By Romanoff's arguments \cite{Romanoff34} (or see
\cite[p.203]{Nathanson96}), we know that
$$
\sum_{2^{n_1}<2^{n_2}\leq N}\prod_{p\mid
2^{n_2-n_1}-1}\bigg(1+\frac{1}{p}\bigg)\ll(\log N)^2.
$$
This concludes the proof of Theorem \ref{sumset}
\end{proof}

\section{Proof of Theorems \ref{sumsetcomplement} and \ref{sumdiffcomplement}}
\setcounter{equation}{0} \setcounter{Thm}{0} \setcounter{Lem}{0}
\setcounter{Cor}{0}

For an integer $a$ and a positive integer $n$, let $a(n)$ denote the
residue class $\{x\in\Z:\, x\equiv a\pmod{n}\}$. For a finite system
$\A=\{a_{s}(n_s)\}_{s=1}^k$, we say $\A$ is a cover of $\Z$ provided
that
$$
\bigcup_{s=1}^ka_s(n_s)=\Z.
$$
Our aim is to find a cover $\{a_{s}(n_s)\}_{s=1}^k$ of $\Z$ and
distinct primes $p_1,p_2,\ldots,p_k$ with $p_s\equiv 3\pmod{4}$ and
$2^{n_s}\equiv1\pmod{p_s}$. With help of the book \cite{BLSTW02},
the following lemma can be directly verified.
\begin{Lem}
\label{coverl} Let
\begin{align*}
\{(a_s',n_s',p_s)\}_{s=1}^{49} =&\{(0,3,7), (1,15,11), (4,15,31),
(7,15,151), (10,15,331),\\ &(13,105,43), (28,105,71), (43,105,127),
(58,105,211),\\ &(73,105,29191),(88,105,86171), (103,315,870031),\\
&(208,315,983431), (313,315,1765891), (2,9,19), (5,27,87211),\\
&(14,81,71119), (41,81,97685839), (68,81,163), (23,135,271),\\
&(50,135,631), (77,135,811), (104,135,23311), (131,135,348031),\\
&(8,99,23), (17,99,67), (26,99,199), (35,99,683), (44,99,5347),\\
&(53,99,599479), (62,99,33057806959),
(71,99,242099935645987),\\
&(80,495,991), (179,495,2971), (278,495,3191),\\
&(377,495,48912491), (476,495,2252127523412251), (89,693,463),\\
&(188,693,5419), (287,693,14323), (386,693,289511839),\\
&(485,693,35532364099), (584,693,2868251407519807),\\
&(683,693,581283643249112959), (98,297,694387),\\
&(197,297,14973866897175265228063698945547), (296,891,1783),\\
&(593,891,1409033313878253109224688819),\\
&(890,891,12430037668834128259094186647)\}
\end{align*}
Then $\A'=\{a_s'(n_s')\}_{s=1}^{49}$ is a cover of $\Z$. And for
$1\leq s\leq 49$, we have $p_s\mid 2^{n_s'}-1$ or $p_s\mid
2^{n_s'}+1$.
\end{Lem}
\begin{Rem}
In \cite{WuSun}, Wu and Sun constructed a cover of $\Z$ with 173 odd
moduli and distinct primitive prime divisors.
\end{Rem}

For $1\leq s\leq 49$, let $n_s=2n_s'$ and let $a_s$ be an integer
such that $a_s\equiv a_s'\pmod{n_s'}$ and $a_s\equiv 1\pmod{2}$. Let
$a_{50}=0$, $n_{50}=2$ and $p_{50}=3$. Then by the Chinese remainder
theorem,
\begin{equation}
\label{cover} \A=\{a_s(n_s)\}_{s=1}^{50}\tag{$*$}
\end{equation} is a cover of $\Z$, and
$2^{n_s}\equiv1\pmod{p_s}$ for $1\leq s\leq 50$. Let
$$
M_1=\prod_{s=1}^{50}p_s
$$
and let $\alpha_1$ be an integer such that
$$\alpha_1\equiv
2^{a_s}\pmod{p_s}$$ for $1\leq s\leq 50$.

Let $x$ be an arbitrary positive integer with $x\equiv
\alpha_1\pmod{M_1}$. Suppose that $x\in\mathcal{S}_3$, i.e.,
$x=p^2+b^2+2^n$ for some $p\in\P$ and $b,n\in\N$. Since $\A$ is a
cover of $\Z$, there exists $1\leq s\leq50$ such that $n\equiv
a_s\pmod{n_s}$. Then
$$
p^2+b^2=x-2^n\equiv \alpha_1-2^{a_s}\equiv 0\pmod{p_s}.
$$
Noting that $p_s\equiv 3\pmod{4}$, $-1$ is a quadratic non-residue
modulo $p_s$. It follows that $$p\equiv b\equiv 0\pmod{p_s}.$$ Since
$p$ is prime, we must have $p=p_s$.

Below we require some additional congruences. Arbitrarily choose
distinct primes $q_1,q_2,\ldots,q_{50}$ such that $(q_s,M_1)=1$ and
$q_s\equiv 7\pmod{8}$ for $1\leq s\leq 50$. Clearly $2$ is a
quadratic residue and $-1$ is a quadratic non-residue modulo $q_s$.
So $-2^n$ is a quadratic non-residue modulo $q_s$ for any $n\geq 0$.
Let
$$
M_2=\prod_{s=1}^{50}q_s,
$$
and let $\alpha_2$ be an integer such that
$$
\alpha_2\equiv p_s^2\pmod{q_s}
$$
for every $1\leq s\leq 50$.

Let $M=M_1M_2$, and let $\alpha$ be an integer such that
$$
\alpha\equiv\alpha_i\pmod{M_i}
$$
for $i=1,2$. Then we have $\{x\in\N:\,
x\equiv\alpha\pmod{M}\}\cap\mathcal{S}_3=\emptyset$. In fact, assume
on the contrary that $x\equiv\alpha\pmod{M}$ and $x=p^2+b^2+2^n$ for
some $p\in\P$ and $b,n\in\N$. Noting that
$x\equiv\alpha_1\pmod{M_1}$, we know $p=p_s$ for some $1\leq s\leq
50$. But since $x\equiv\alpha_2\pmod{M_2}$,
$$
x-p_s^2-2^n\equiv\alpha_2-p_s^2-2^n\equiv-2^n\pmod{q_s}.
$$
So $x-p_s^2-2^n$ is a quadratic non-residue modulo $q_s$, which
leads to an evident contradiction since $x-p_s^2-2^n=b^2$. This
concludes the proof of Theorem \ref{sumsetcomplement}.\qed

Now let us turn to the proof of Theorem \ref{sumdiffcomplement}. We
still use the cover $\A$ in (\ref{cover}). Now suppose that $x\equiv
-\alpha_1\pmod{M_1}$ and there exist $p\in\P$ and $b,n\in\Z$ such
that $x=p^2+b^2-2^n$. Then $n\equiv a_s\pmod{n_s}$ for some $1\leq
s\leq 50$, and
$$
x+2^n\equiv -\alpha_1+2^{a_s}\equiv 0\pmod{p_s}.
$$
It follows that $p=p_s$. The main difficult is to find the
additional congruences.
\begin{Lem}
\label{addcong} Let
\begin{align*}
\{(c_s,r_s)\}_{s=1}^{50} =&\{(505,47\times 178481),(5519,601\times
1801), (366,2731\times 8191),\\ &(1303,73\times 262657),
(5149,233\times 2089), (5938,43691\times 131071),\\
&(182725,223\times 616318177), (12153,174763\times 524287),\\
&(148671,13367\times 164511353), (490297,431\times 2099863),\\
&(115115,2351\times 13264529), (2370639,6361\times 20394401),\\
&(37,5\times 17\times 257), (5615,13\times 37\times 109),
(146,89\times 397\times 2113),\\
&(637,97\times 241\times 673), (6393,103\times 2143\times 11119),\\
&(13847,53\times 157\times 1613), (1799,29\times 113\times
15790321),\\
&(335,59\times 1103\times 3033169), (451,337\times 92737\times
649657),\\
&(1479,641\times 65537\times 6700417), (40655,137\times 953\times
26317),\\
&(13353,228479\times 48544121\times 212885833),\\ &(23775,439\times
2298041\times 9361973132609), (10334,229\times 457\times 525313),\\
&(65971,2687\times 202029703\times 1113491139767), (5893,41\times
61681\times 4278255361),\\ &(1344867,911\times 112901153\times
23140471537), (560826,277\times 1013\times 30269),\\
&(406789,283\times 4513\times 165768537521),\\ &(415099,191\times
420778751\times 30327152671),\\
&(61153,101\times 4051\times 8101), (1261375,307\times 2857\times
6529),\\
&(1324442,107\times 69431\times 28059810762433),\\
&(1663519,321679\times 26295457\times 319020217),\\
&(2094571,3391\times 23279\times 65993), (1032375,571\times
32377\times 1212847),\\
&(6321391,14951\times 4036961\times 2646507710984041),\\
&(19031871,937\times 6553\times 7830118297), (7918330,2833\times
37171\times 179951),\\
&(2286429,61\times 1321\times 4562284561), (2227201,5581\times
8681\times 49477),\\
&(207773684,131\times 409891\times 7623851), (2526613,281\times
122921\times 7416361),\\
&(5596695,433\times 577\times 38737), (25234915,593\times 1777\times
25781083),\\
&(7950774,251\times 100801\times 10567201), (10130779,313\times
21841\times 121369),\\
&(14272093,1429\times 3361\times 14449)\}.
\end{align*}
Then for every $1\leq s\leq 50$ and $n\in\N$, $c_s+2^n$ is a
quadratic non-residue modulo $r_s$.
\end{Lem}
Lemma \ref{addcong} can be checked via a direct computation. In
fact, we only need to consider those $c_s+2^n$ modulo $r_s$ for
$0\leq n<\ord_2(r_s)$, where $\ord_2(r)$ denotes the least positive
integer such that $2^{\ord_2(r)}\equiv 1\pmod{r}$.

Notice that $(r_i,M_1)=1$ and $(r_i,r_j)=1$ for any distinct $i,j$.
Let
$$
M_3=\prod_{s=1}^{50}r_s
$$
and let $\alpha_3$ be an integer such that
$$
\alpha_3\equiv c_s+p_s^2\pmod{r_s}.
$$
Let $M'=M_1M_3$ and let $\alpha'$ be an integer satisfying
$$
\alpha\equiv-\alpha_1\pmod{M_1}\qquad\text{and}\qquad
\alpha\equiv\alpha_3\pmod{M_3}.
$$
For any $x\equiv\alpha'\pmod{M'}$, assume on the contrary that
$x=p_s^2+b^2-2^n$ for some $1\leq s\leq 50$ and $b,n\in\N$. Then
$$
b^2=x+2^n-p_s^2\equiv\alpha_3+2^n-p_s^2\equiv c_s+2^n\pmod{r_s}.
$$
This is impossible since $c_s+2^n$ is a quadratic non-residue modulo
$r_s$. Hence the reside class $\{x\in\N:\,
x\equiv\alpha'\pmod{M'}\}$ contains no integer of the form
$p^2+b^2-2^n$.\qed
\begin{Rem}
Observe that the moduli appear in Lemma \ref{addcong} are all
composite. So we have the following problem.
\begin{Prob}
Does there exist infinitely many primes $p$ such that the set
$$
\{1\leq c\leq p:\, c+2^n\text{ is a quadratic non-residue modulo
}p\text{ for every }n\in\N\}
$$
is non-empty?
\end{Prob}
In fact, we don't know any such prime $p$. For example, let
$p=2^{19}-1$, then $27006+2^1,27006+2^2,\ldots,27006+2^{18}$ are all
quadratic non-residues modulo $p$, but $27006+2^{19}$ is a quadratic
residue modulo $p$.

\end{Rem}

\section{The integers of the form $b_1^2+b_2^2+2^{n^2}$}
\setcounter{equation}{0} \setcounter{Thm}{0} \setcounter{Lem}{0}
\setcounter{Cor}{0}
\begin{proof}[Proof of the first assertion of Theorem \ref{sumtwosqaures}]
Let
$$
\mathcal{Q}=\{x\in\N:\, x\text{ has no prime factor of the form
}4k+3\}.
$$
Clearly
$$
\mathcal{Q}\subseteq\{b_1^2+b_2^2:\, b_1,b_2\in\N\}.
$$
We only need to prove that the set
$$
\{x+2^{n^2}:\,x\in\mathcal{Q},\ n\in\N\}
$$
has a positive lower density. As an application of half dimensional
sieve method \cite{Iwaniec76}, we know that
$$
|\{x\in\mathcal{Q}:\,x\leq N\}|\gg\frac{N}{\sqrt{\log N}}.
$$
And with help of Selberg's sieve method, it is not difficult to see
that
$$
|\{(x_1,x_2):\, x_1=x_2+m,\ x_i\in\mathcal{Q},\
x_i\leq N\}|\ll\frac{N}{\log N}\prod_{\substack{p\mid m\\
p\equiv 3\pmod{4}}}\bigg(1+\frac{1}{p}\bigg)
$$
for every positive integer $M$. By the Cauchy-Schwarz inequality,
\begin{align*}
&|\{x+2^{n^2}:\,x+2^{n^2}\leq N,\ x\in\mathcal{Q},\
n\in\N\}|\\
\geq&\frac{|(x,n):\,x+2^{n^2}\leq N,\ x\in\mathcal{Q},\
n\in\N\}|^2}{|\{(x_1,x_2,n_1,n_2):\,x_1+2^{n_1^2}=x_2+2^{n_2^2}\leq
N,\ x_i\in\mathcal{Q},\ n_i\in\N\}|}.
\end{align*}
So it suffices to show that
$$
|\{(x_1,x_2,n_1,n_2):\,x_1+2^{n_1^2}=x_2+2^{n_2^2}\leq N,\
x_i\in\mathcal{Q},\ n_i\in\N\}|\ll N.
$$
Now
\begin{align*}
&|\{(x_1,x_2,n_1,n_2):\,x_1+2^{n_1^2}=x_2+2^{n_2^2}\leq N,\
x_i\in\mathcal{Q},\ n_i\in\N\}|\\
\leq&|\{(x_1,n_1):\,x_1+2^{n_1^2}\leq N,\ x_1\in\mathcal{Q},\
n_1\in\N\}|\\
&+2\sum_{\substack{0\leq n_1<n_2\leq\sqrt{\log N/\log 2}}}
|\{(x_1,x_2):\,x_1-x_2=2^{n_2^2}-2^{n_1^2},\
x_i\in\mathcal{Q}\cap[1,N]\}|\\
\ll&N+\frac{N}{\log N}\sum_{\substack{0\leq n_1<n_2\leq\sqrt{\log
N/\log 2}}}\prod_{\substack{p\mid 2^{n_2^2}-2^{n_1^2}\\ p\equiv
3\pmod{4}}}\bigg(1+\frac{1}{p}\bigg).
\end{align*}
Obviously
\begin{align*}
\sum_{\substack{0\leq n_1<n_2\leq\sqrt{\log N/\log
2}}}\prod_{\substack{p\mid 2^{n_2^2}-2^{n_1^2}\\ p\equiv
3\pmod{4}}}\bigg(1+\frac{1}{p}\bigg)\leq&\sum_{\substack{0\leq
n_1<n_2\leq\sqrt{\log N/\log 2}}}\prod_{\substack{p\mid
2^{n_2^2}-2^{n_1^2}}}\bigg(1+\frac{1}{p}\bigg)\\
\leq&\sum_{d}\frac{1}{d}\sum_{\substack{0\leq n_1<n_2\leq\sqrt{\log
N/\log 2}\\ n_2^2\equiv n_1^2\pmod{\ord_2(d)}}}1\\
\end{align*}
Suppose that $p$ is prime, $\beta\geq 1$ and $1\leq a\leq p^\beta$.
Then we have
$$
|\{1\leq x\leq p^\beta:\, x^2\equiv a\pmod{p^\beta}\}|\leq
2p^{\frac{\nu_p(a)}{2}}
$$
since the multiplicative group modulo $p^\beta$ is cyclic, where
$\nu_p(a)$ denotes the greatest integer such that $p^{\nu_p(a)}\mid
a$. Thus
\begin{align*}
\sum_{\substack{0\leq n_1<n_2\leq\sqrt{\log N/\log 2}\\ n_2^2\equiv
n_1^2\pmod{\ord_2(d)}}}1\ll\sqrt{\frac{\log N}{\log
2}}\bigg(\sqrt{\frac{\log N}{\log
2}}\cdot\frac{2^{\omega(\ord_2(d))}\sqrt{\ord_2(d)}}{\ord_2(d)}+1\bigg),
\end{align*}
where $\omega(r)$ denotes the number of distinct prime factors of
$r$. We only need to prove that
$$
\sum_{d}\frac{2^{\omega(\ord_2(d))}}{d\sqrt{\ord_2(d)}}
$$
converges. Define
$$
E(x)=\sum_{k\leq x}\sum_{\ord_2(d)=k}\frac{1}{d}.
$$
Romanoff had show that $E(x)\ll\log x$ (cf. \cite[pp.
200-2001]{Nathanson96}). So
\begin{align*}
\sum_{d}\frac{2^{\omega(\ord_2(d))}}{d\sqrt{\ord_2(d)}}=&
\sum_{k}\frac{2^{\omega(k)}}{\sqrt{k}}\sum_{\ord_2(d)=k}\frac{1}{d}\\
\ll&\int_{1}^\infty x^{-\frac{1}{3}}d(E(x))\\
=&x^{-\frac{1}{3}}E(x)\bigg|_{1}^\infty+\frac{1}{3}\int_{1}^\infty
x^{-\frac{4}{3}}E(x)dx=O(1).
\end{align*}
Our proof is complete.
\end{proof}
Let $\A=\{a_s(n_s)\}_{s=1}^{50}$ be the cover in (\ref{cover}), and
let $p_1,\ldots,p_{50}$ be the corresponding primes with $p_s\equiv
3\pmod{4}$ and $p_s\mid 2^{n_s}-1$. Since every $p_s$ has at least
one quadratic non-residue modulo $p_s$, the second assertion of
Theorem \ref{sumtwosqaures} is an immediate consequence of the
following stronger result.
\begin{Thm}
\label{twosquare} Let $\mathcal{N}$ be a set of non-negative
integers. Suppose that for every $1\leq s\leq 50$, there exists
$1\leq e_s\leq p_s$ such that
$$
|\mathcal{N}\cap\{x\in\N:\, x\equiv e_s\pmod{p_s}\}|<+\infty.
$$
Then there exists a residue class with odd modulo, which contains no
integer of the form $b_1^2+b_2^2+2^n$ with $n\in\mathcal{N}$.
\end{Thm}
\begin{proof}
Let $n_s^*=\ord_2(p_s)$. Noting that $n_s^*\mid n_s$ and
$(n_s,p_s)=1$, for every $1\leq s\leq 50$, let $a_s^*$ be an integer
such that
$$
a_s^*\equiv a_s\pmod{n_s^*}
$$
and
$$
a_s^*\equiv e_s\pmod{p_s}.
$$
Clearly $\A^*=\{a_s^*(n_s^*)\}_{s=1}^{50}$ is also a
cover of $\Z$.

Let
$$
\mathscr{H}_s=\mathcal{N}\cap\{x\in\N:\, x\equiv
e_s\pmod{p_s}\}=\{h_{s,1},h_{s,2},\ldots,h_{s,|\mathscr{H}_s|}\}.
$$
for $1\leq s\leq 50$. Choose distinct
$|\mathscr{H}_1|+|\mathscr{H}_2|+\cdots+|\mathscr{H}_{50}|$ primes
$$
q_{1,1},\ \ldots,\ q_{1,|\mathscr{H}_1|},\ q_{2,1},\ \ldots,\
q_{2,|\mathscr{H}_2|},\ \ldots,\ q_{50,1},\ \ldots,\
q_{50,|\mathscr{H}_{50}|}
$$
satisfying that
$$
q_{s,t}\equiv 3\pmod{4}
$$
and
$$
q_{s,t}\not\in\{p_1,p_2,\ldots,p_{50}\}
$$
for every $1\leq s\leq 50$ and $1\leq t\leq |\mathscr{H}_{s}|$.

Let
$$
M^*=\bigg(\prod_{1\leq s\leq 50}p_s\cdot\prod_{\substack{1\leq s\leq 50\\
1\leq t\leq |\mathscr{H}_{s}|}}q_{s,t}\bigg)^2
$$
and let $\alpha^*$ be an integer such that
$$
\alpha^*\equiv 2^{a_s^*}\pmod{p_s^2}$$ and
$$
\alpha^*\equiv 2^{h_{s,t}}+q_{s,t}\pmod{q_{s,t}^2}$$ for every
$s,t$.

We claim that for any $x\equiv\alpha^*\pmod{M^*}$, $x$ is not of the
form $b_1^2+b_2^2+2^{n}$ with $n\in\mathcal{N}$. Assume on the
contrary that $x\equiv\alpha^*\pmod{M^*}$ and $x=b_1^2+b_2^2+2^{n}$
with $n\in\mathcal{N}$. Since $\A^*=\{a_s^*(n_s^*)\}_{s=1}^{50}$ is
a cover, similarly as the arguments in the proof of Theorems
\ref{sumsetcomplement} and \ref{sumdiffcomplement}, we know that
$$
b_1\equiv b_2\equiv 0\pmod{p_s}
$$
for some $1\leq s\leq 50$. It follows that
$$
x-2^{n}\equiv 2^{a_s^*}-2^{n}\equiv 0\pmod{p_s^2},
$$
that is, $n\equiv a_s^*\pmod{\ord_{2}(p_s^2)}$. It is not difficult
to check that
$$
2^{n_s^*}=2^{\ord_2(p_s)}\not\equiv 1\pmod{p_s^2}.
$$
(In fact, the only known primes $p$ with $2^{p-1}\equiv1\pmod{p^2}$
are 1093 and 3511.) And
$$
2^{n_s^*p}=\sum_{k=0}^{p}\binom{p}{k}(2^{n_s^*}-1)^k\equiv
1\pmod{p_s^2}.
$$
So we must have $\ord_2(p_s^2)=n_s^*p$. Consequently
$$
n\equiv a_s^*\equiv e_s\pmod{p_s}.
$$
Since $n\in \mathcal{N}$, we have $n\in\mathscr{H}_{s}$ and there
exists $1\leq t\leq |\mathscr{H}_{s}|$ such that $n=h_{s,t}$. It
follows that
$$
b_1^2+b_2^2=x-2^{n}\equiv\alpha^*-2^{h_{s,t}}\equiv
q_{s,t}\equiv0\pmod{q_{s,t}^2}.
$$
Recalling that $q_{s,t}\equiv 3\pmod{4}$, $q_{s,t}\mid b_1^2+b_2^2$
implies that
$$
b_1\equiv b_2\pmod{q_{s,t}}.
$$
But it is impossible since
$$
b_1^2+b_2^2\equiv q_{s,t}\not\equiv 0\pmod{q_{s,t}^2}.
$$
\end{proof}

\begin{Cor}
There exists a positive integer $m$ such that the set
$$
\{x\in\N:\, x\text{ is even and }x\text{ is not of the form
}b_1^2+b_2^2+2^{mn}\}
$$
contains an infinite arithmetic progression.
\end{Cor}
\begin{proof}
Let $m=p_1p_2\ldots p_{50}$ where $p_1,p_2,\ldots,p_{50}$ are the
primes in Lemma \ref{coverl}. Thus substituting
$\mathcal{N}=\{x\in\N:\, x\equiv0\pmod{m}\}$ and $e_s=1$ in Theorem
\ref{twosquare}, we are done.
\end{proof}
\begin{Prob} Does there exists a residue class with odd modulo, which contains no
integer of the form $b_1^2+b_2^2+2^n$ with $b_1,b_2,n\in\N$?
\end{Prob}

\begin{Ack}
We are grateful to Professors Hongze Li and Zhi-Wei Sun for their
helpful suggestions.
\end{Ack}

\end{document}